\documentclass[reqno]{amsproc}

\usepackage{amstext,amsmath,amssymb,amsfonts}
\usepackage[latin1]{inputenc}
\usepackage{epsfig}
\usepackage{color}
\usepackage{tikz}
\usepackage[all,arc,dvips,arrow,curve,cmactex]{xy}

\usepackage{latexsym}
\usepackage{mathtools}
\DeclarePairedDelimiter\floor{\lfloor}{\rfloor}

\newtheorem{lemma}{Lemma}

\newtheorem{definition}{Definition}
\newtheorem{theorem}{Theorem}

\newtheorem{proposition}{Proposition}
\newtheorem{example}{Example}

\newtheorem{conjecture}{Conjecture}

\newcommand{\bea}{\begin{eqnarray}}
\newcommand{\eea}{\end{eqnarray}}
\newcommand{\beq}{\begin{equation}}
\newcommand{\eeq}{\end{equation}}
\newcommand{\enn}{\nonumber \end{equation}}

 \newcommand{\rk}{{\rm r}\,}

 \newcommand{\cF}{\mathcal{F}}

 \newcommand{\cB}{\mathcal{B}}

\title[Canonical $\Delta$-matroid]{CANONICAL BINARY $\Delta$-MATROIDS}

\author{R\'emi Cocou Avohou}
\address[R.C.A.]{
Humboldt-Universit{\"a}t zu Berlin, Institut f{\"u}r Mathematik und Institut f{\"u}r Physik, Rudower Chaussee
25, 12489 Berlin, Germany, \& ICMPA-UNESCO Chair, 072BP50, Cotonou,
Rep. of Benin, \& Ecole Normale Superieure, B.P 72, Natitingou, Benin}
\email{avohoure@hu-berlin.de}

\author{Brigitte Servatius}
\address[B.S.]{Mathematical Sciences, Worcester Polytechnic Institute, Worcester MA 01609-
2280}
\email{bservat@wpi.edu }

\author{Herman Servatius}
\address[B.S.]{Mathematical Sciences, Worcester Polytechnic Institute, Worcester MA 01609-
2280}
\email{hservat@wpi.edu }
\begin{document}

\maketitle

\begin{abstract}
The handle slide operation, originally defined for ribbon graphs, was extended to 
delta-matroids by  I. Moffatt and E. Mphako-Bandab, who 
show that, using a delta-matroid analogue of handle slides, every binary delta-matroid in which the empty set is feasible can be written in a canonical form analogous to the canonical form for one-vertex maps on a surface.
We provide a canonical form for binary delta-matroids without restriction on the feasibility of the empty set.

%
%
\end{abstract}
\tableofcontents
\section{Introduction}
Whitney introduced the concept of matroid in 1935, while working on abstract properties of linear dependence, and by now there are many  excellent books on matroid theory, see for example \cite{MR0427112,MR0276117,MR1027839,MR3154793, Oxley11}. Matroid theory has an interesting feature in that there are many different but equivalent ways of defining a matroid: in terms of independent sets, definition in terms of circuits, in terms of bases, in terms of spanning sets, to name just a few.

In terms of bases, a matroid is defined as follows. A matroid $M$ is a pair $(E,\cB)$, where $E$ is a nonempty finite set and $\cB$ is a nonempty collection of bases, which are subsets of $E$,  satisfying the following requirement known as the basis exchange axiom (EA).

(EA) if $B_1$ and $B_2$ are bases, and $b_1$ is any element of $B_1\setminus B_2$, then there is  an element $b_2$ of $B_2\setminus B_1$ such $(B_1\setminus \{b_1\})\cup \{b_2\}=B_1\Delta\{b_1,b_2\}$ is also a base. 

Subsets of bases are referred to as \emph{ independent sets}, sets that are not independent are called \emph{dependent sets}, minimal dependent sets are called \emph{cycles}, and sets containing a basis are referred to as \emph{spanning}.

A matroid is called \emph{representable} over a field $F$ if its elements can be injectively mapped into a set of vectors of a vector space over $F$ such that independent sets get mapped to linearly independent sets of vectors. A matroid is \emph{binary} if it can be represented by a vector space over $GF(2)$, the finite field on two elements.  A regular matroid is representable  over every field.

By substituting the symmetric difference for the set difference in the Exchange Axiom (EA), we obtain the Symmetric Exchange Axiom (SEA), which is used by Bouch\`{e}t \cite{MR904585,MR1020647,MR1083464} to define $\Delta$-matroids  extending the concept of matroid \cite{Oxley11}. A $\Delta$-matroid $D$ is a finite set $E$ and a collection $\cF$ of subsets of $E$ called \emph{feasible sets} satisfying the condition that

(SEA) if $F_1$ and $F_2$ are in $\cF$ and $x \in F_1\Delta F_2$ then there exists a $y \in F_2\Delta F_1$ such that $F_1\Delta \{x, y\} \in \cF$. Note  that $x=y$ is allowed.

The  notion of binary matroid is naturally extended to $\Delta$-matroids by A. Bouch\`{e}t in  \cite{MR904585,MR1020647,MR1083464,MR4193372,AS2019} giving birth to binary $\Delta$-matroids.  Furthermore, an interesting operation known in the class of one vertex ribbon graphs, also known as rosettes or disc-band surfaces, found an extension to $\Delta$-matroids \cite{MR3546900}.
And, the class of binary 
$\Delta$-matroids is proved to be the only class of $\Delta$-matroids closed under handle slides \cite{MR4207277}.
A handle slide in the context of ribbon graphs is the movement of the end of one edge over an edge adjacent to it in cyclic order at a vertex. Under this operation, every rosette has the canonical form $B_{i,j,k}$, which consists of $i$ orientable loops not interlacing any other loop, $j$ pairs of interlaced orientable loops, and $k$ non-orientable loops not interlacing any other loop. This is essentially the classification of surfaces with boundary up to homeomorphism: $j$ corresponds to the number of tori, $k$ the number of real projective planes, and $i+1$ the number of holes in the surface.

In \cite{MR3546900}, a more algebraic formulation of the handle slides operation was introduced, and an analogue of canonical rosette was found for the class of binary $\Delta$-matroids where the empty set is a feasible:

\begin{theorem}\label{theo:binarydeltaempty}
Let $D = (E, \cF)$ be a binary $\Delta$-matroid in which the empty set is feasible. Then, for some $i$, $j$, $k$, there is a sequence of handle slides taking $D$ to $D_{i,j,0}$ if $D$ is even, or $D_{i,0,k}$, with
$k\neq 0$, if $D$ is odd. Furthermore, if some sequences of handle slides take $D$ to $D_{i,j,k}$ and to $D_{p,q,r}$ then $i = p$, and so $D$ is taken to a unique form $D_{i,j,0}$ or $D_{i,0,k}$ by handle slides.
\end{theorem}

In this theorem each $D_{i,j,k}$ represents the $\Delta$-matroid of $B_{i,j,k}$ in which the ground set is the edges of $B_{i,j,k}$ and the feasible set is the collection of its spanning quasi-trees.

The same paper \cite{MR3546900} asked for a canonical form for a binary $\Delta$-matroid and  conjecture such a form. This work gives the proof of a refined version. The section \ref{sect:prelimi} that follows provides a quick review of some fundamentals in matroids and $\Delta$-matroids. Section \ref{sect:conproof} deals with a set of lemmas as well as the proof of the conjecture that is reformulated in Theorem \ref{theo:deltamatcanoni}.
\section{Preliminaries}
\label{sect:prelimi}

%
Let $D=(E, \cF)$ be a delta-matroid.
There are two matroids found as substructures of $D$, namely the \emph{upper matroid} whose bases are the feasible sets of $D$ of largest cardinality, and the \emph{lower matroid}, whose bases are the smallest feasible sets of $D$.

Let $M=(E,\cB)$ be a matroid defined by its bases set $\cB$. Subsets of bases are called independent sets, sets that are not independent are dependent sets, minimal dependent sets are called cycles, and sets containing a basis are referred to as spanning sets. We have the following theorem from \cite{AS2021,AS2019}.

\begin{theorem}\label{the:spanind}
If $D$ is a $\Delta$-matroid and $F$ is a feasible set of $D$, then $F$ is spanning
in the lower matroid and independent in the upper matroid.
\end{theorem}

To each matroid $M = (E, \cB) $ there is a \emph{dual} matroid $M^* = (E, \cB^*)$ with $\cB^* = \{E\setminus B | B \in \cB$, and to each $\Delta$-matroid $D = (E, \cF) $ there is a  \emph{dual} $\Delta$-matroid $D^* = (E, \cF^*)$ with 
$\cF^* = \{E\setminus F | F \in \cF $. A \emph{loop} of $D$ is an element of $E$
which is not contained in any feasible set of $D$, while a \emph{coloop} of $D$ is not contained in any feasible set of $D^*$. Note that the upper (lower) matroid of $D^*$ is the dual of the lower (upper) matroid of $D$.

\begin{definition}[Elementary minors] \label{def:minor}
Let $D=(E, \cF)$ be a delta-matroid. The elementary minors of $D$ at $e\in E$, are the delta-matroids $D-e$ and $D/e$ defined by:
$$D-e=\Big(E-e, \big\{F | F\subseteq E-e, F\in \cF\big\}\Big),$$
if $e$ is not a coloop, and 
$$D/e=\Big(E-e, \big\{F | F\subseteq E-e, F\cup e\in \cF\big\}\Big),$$
if $e$ is not a loop. In case $e$ is a loop or a coloop, we set $D/e=D-e$.
The delta-matroid $D-e$ is called the deletion of $D$ along $e$, and $D/e$ the contraction of $D$ along $e$.
\end{definition}

A $\Delta$-matroid obtained from a $\Delta$-matroid $D$ by a ( possibly empty) sequence of deletions and contractions is called a \emph{minor} of $D$.

Let $A=(a_{vw}: v, w\in E)$ be a symmetric binary matrix and $A[W]=(a_{vw}: v, w\in W)$ for $W\subseteq E$. 
Assuming that $A[\emptyset]$ has an inverse, $D(A)=(E, \{W: A[W] \text{ has an inverse}\})$ is a $\Delta$-matroid.

\begin{definition}[Twist]\label{def:twist}
Let $D=(E, \cF)$ be a set system. For $A\subseteq E$, the twist of $D$ with respect to $A$, denoted by $D\star A$, is given by $(E, \{A\Delta X | X\in \cF\})$.
\end{definition}

Note that for  the dual  $D^*$  of $D$ we have  $ D^* = D\star E$.

\begin{definition}[Binary delta-matroid \cite{MR1020647, MR1083464}]\label{def:binarydeltamat}
A  delta-matroid $D=D(E, \cF)$ is said to be binary if there exists $F\in \cF$ and a symmetric binary matrix $A$ such that $D=D(A)\star F$.
\end{definition}

\begin{definition}[Handle slides \cite{MR3546900}]
\label{def:hs}
Let $D=(E, \cF)$ be a set system, and $a, b\in E$ with $a\neq b$. We define $D_{ab}$ to be the set system $(E, \cF_{ab})$ where
$$\cF_{ab}=\cF\Delta \big\{X\cup a | X\cup b\in\cF \text{ and } X\subseteq E-\{a, b\}\big\}.$$
We call the move taking $D$ to $D_{ab}$ a handle slide taking $a$ over $b$.
\end{definition}

Performing the handle slide operation on a given $\Delta$-matroid does not necessarily result in another $\Delta$-matroid but the result is stable only for the class of binary $\Delta$-matroids \cite{MR4207277}. There is a sequence of handle slides that can send any binary $\Delta$-matroid, in which the empty set is a feasible, to a particular $\Delta$-matroid called canonical binary $\Delta$-matroid.

The canonical binary $\Delta$-matroid denoted by $D_{i,j,k}$ is a binary $\Delta$-matroid in which the empty set is a feasible and that arises as the direct sum of $i$ copies of $({e}, \{\emptyset\})$, $j$ copies of $(\{e, f\}, \{\emptyset, \{e, f\}\})$, and $k$ copies of $(\{e\}, \{\emptyset, \{e\}\})$. It is worth noting that the sum is performed on isomorphic copies of these $\Delta$-matroids with mutually disjoint ground sets. We denote by $D_{i,j,k,l}$ the $\Delta$-matroid
consisting of the direct sum of $D_{i,j,k}$ with $l$ copies of the $\Delta$-matroids isomorphic to $({e}, \{\{e\}\})$.

From this the conjecture is as follows:
\begin{conjecture}
For each binary $\Delta$-matroid $D$, there is a sequence of handle slides taking $D$ to some $D_{i,j,k,l}$ where $i$ is the size of the ground set minus the size of a largest feasible set, $l$ is the size of a smallest feasible set, $2j + k$ is difference in the sizes of a largest and a smallest feasible
set. Moreover, $k = 0$ if and only if D is even, and if $D$ is odd then every value of $j$ from $0$ to $\floor*{\frac{w}{2}}$,
where $w$ is the difference between the sizes of a largest and a smallest feasible set, can be attained.
\end{conjecture}

\section{Canonical binary matroids}
\label{sect:conproof}
Since a $Delta$-matroid for which upper and lower matroid are the same is just a matroid, we start by looking for a normal form for binary matroids.

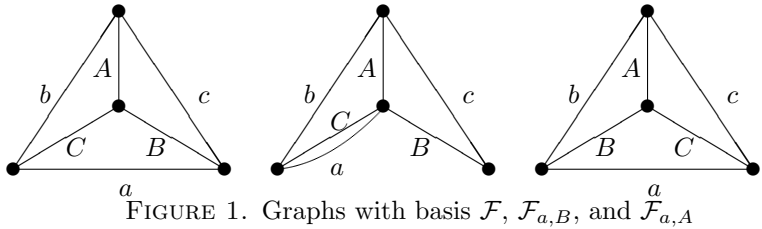
\begin{figure}[htb]
\begin{picture}(300,70)(0,5)
 \put(0,0){\circle*{5}}
 \put(80,0){\circle*{5}}
 \put(0,0){\line(1,0){80}}
 \put(0,0){\line(2,3){40}}
 \put(0,0){\line(5,3){40}}
 \put(40,24){\circle*{5}}
 \put(40,24){\line(0,1){36}}
 \put(40,60){\line(2,-3){40}}
 \put(40,24){\line(5,-3){40}}
 \put(40,60){\circle*{5}}
 \put(40,-10){\mbox{$a$}}
 \put(10,25){\mbox{$b$}}
 \put(70,25){\mbox{$c$}}
 \put(20,5){\mbox{$C$}}
 \put(50,5){\mbox{$B$}}
 \put(30,35){\mbox{$A$}}

\bezier{100}(140,23)(120,2)(100,0)
\put(100,0){\circle*{5}}
 \put(180,0){\circle*{5}}
 \put(100,0){\line(2,3){40}}
 \put(100,0){\line(5,3){40}}
 \put(140,24){\circle*{5}}
 \put(140,24){\line(0,1){36}}
 \put(140,60){\line(2,-3){40}}
 \put(140,24){\line(5,-3){40}}
 \put(140,60){\circle*{5}}
 \put(120,-2){\mbox{$a$}}
 \put(110,25){\mbox{$b$}}
 \put(170,25){\mbox{$c$}}
 \put(120,15){\mbox{$C$}}
 \put(150,5){\mbox{$B$}}
 \put(130,35){\mbox{$A$}}

 \put(200,0){\circle*{5}}
 \put(280,0){\circle*{5}}
 \put(200,0){\line(1,0){80}}
 \put(200,0){\line(2,3){40}}
 \put(200,0){\line(5,3){40}}
 \put(240,24){\circle*{5}}
 \put(240,24){\line(0,1){36}}
 \put(240,60){\line(2,-3){40}}
 \put(240,24){\line(5,-3){40}}
 \put(240,60){\circle*{5}}
 \put(240,-10){\mbox{$a$}}
 \put(210,25){\mbox{$b$}}
 \put(270,25){\mbox{$c$}}
 \put(220,5){\mbox{$B$}}
 \put(250,5){\mbox{$C$}}
 \put(230,35){\mbox{$A$}}

\end{picture}
\caption{Graphs with basis $\cF$, $\cF _{a,B}$, and $\cF _{a,A}$  }
\label{tetfig}
\end{figure}

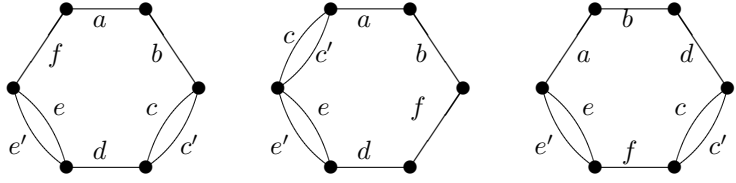
\begin{figure}
\begin{picture}(270,60)(0,0)
 \put(20,0){\circle*{5}}
 \put(50,0){\circle*{5}}
  \put(20,60){\circle*{5}}
 \put(50,60){\circle*{5}}
  \put(0,30){\circle*{5}}
 \put(70,30){\circle*{5}}
\put(20,0){\line(1,0){30}}
\put(20,60){\line(1,0){30}}
\put(0,30){\line(2,3){20}}
\put(70,30){\line(-2,3){20}}
 \put(30,53){\mbox{$a$}}
 \put(30,3){\mbox{$d$}}
 \put(13,40){\mbox{$f$}}
 \put(52,40){\mbox{$b$}}
 \put(15,20){\mbox{$e$}}
 \put(50,20){\mbox{$c$}}
 \put(-2,5){\mbox{$e'$}}
 \put(63,5){\mbox{$c'$}}
 \bezier{100}(0,30)(15,20)(20,0)
 \bezier{100}(0,30)(5,10)(20,0)
 \bezier{100}(70,30)(55,20)(50,0)
 \bezier{100}(70,30)(65,10)(50,0)

\put(220,0){\circle*{5}}
 \put(250,0){\circle*{5}}
  \put(220,60){\circle*{5}}
 \put(250,60){\circle*{5}}
  \put(200,30){\circle*{5}}
 \put(270,30){\circle*{5}}
\put(220,0){\line(1,0){30}}
\put(220,60){\line(1,0){30}}
\put(200,30){\line(2,3){20}}
\put(270,30){\line(-2,3){20}}
 \put(230,53){\mbox{$b$}}
 \put(230,3){\mbox{$f$}}
 \put(213,40){\mbox{$a$}}
 \put(252,40){\mbox{$d$}}
 \put(215,20){\mbox{$e$}}
 \put(250,20){\mbox{$c$}}
 \put(197,5){\mbox{$e'$}}
 \put(263,5){\mbox{$c'$}}
 \bezier{100}(200,30)(215,20)(220,0)
 \bezier{100}(200,30)(205,10)(220,0)
 \bezier{100}(270,30)(255,20)(250,0)
 \bezier{100}(270,30)(265,10)(250,0)
 
 \put(120,0){\circle*{5}}
 \put(150,0){\circle*{5}}
  \put(120,60){\circle*{5}}
 \put(150,60){\circle*{5}}
  \put(100,30){\circle*{5}}
 \put(170,30){\circle*{5}}
\put(120,0){\line(1,0){30}}
\put(120,60){\line(1,0){30}}
\put(170,30){\line(-2,3){20}}
\put(150,0){\line(2,3){20}}
 \put(130,53){\mbox{$a$}}
 \put(130,3){\mbox{$d$}}
 \put(102,47){\mbox{$c$}}
 \put(114,40){\mbox{$c'$}}
 \put(152,40){\mbox{$b$}}
 \put(115,20){\mbox{$e$}}
 \put(150,20){\mbox{$f$}}
 \put(97,5){\mbox{$e'$}}
  \bezier{100}(120,60)(105,50)(100,30)
 \bezier{100}(120,60)(115,40)(100,30)
 \bezier{100}(100,30)(115,20)(120,0)
 \bezier{100}(100,30)(105,10)(120,0)

\end{picture}
\caption{Graphs with identical bases set $\cF$   }
 \label{idfig}
\end{figure}

The term handle slide was originally used to describe the move on ribbon graphs which slides the end of one edge over an edge adjacent to it in the cyclic order at a vertex. In case of 2-connected graphic matroids one can still interpret  Definition~\ref{def:hs} as a sliding move on incident edges, see Figure~\ref{tetfig}, while on non-incident edges the move results in merely a relabelling of the edges. However, if the matroid information is not enough to draw the graph uniquely, for example in the graph in Figure~\ref{idfig}, where sliding $a$ over an incident edge has essentially the same effect as sliding it over a non-incident edge.

In Figures~\ref{tetfig} and~\ref{idfig} we consider the familiar graphic matroid on the edge set of a (connected) graph, whose bases are the spanning trees, which we may also consider the feasible sets of a $\Delta$-matroid. It is well known that the graphic matroid is regular, in particular it is binary, hence an example of a binary $\Delta$-matroid.

Let $M=(E, \cB)$ be a matroid defined by its base set, $F, F'\in \cB$. Then there is a matching between elements of $F\setminus F'$ and $F'\setminus F$ such that $F\Delta \{ x,x'\} \in \cB$  for every matched pair $\{x,x'\}$, which follows directly from EA.
We conclude that every basis  $F'\in \cB$  of $M$ may be written as
$F'=F\Delta \{x_1, y_1\}\Delta \{x_2, y_2\}\Delta\cdots\Delta \{x_P, y_P\}$, with $p\leq |F|$ and such that $x_i\in F$, $y_i\notin F$ $\forall i$.
We say that $p$ is the length of $F'$ in the matroid $M$, and that $\cB$ is described by $F$.

Note that the maximal length of a basis described by $F$ depends on $F$. For example, if $M$ is the cycle matroid of the graph in Figure~\ref{tetfig}, $K_4$, the maximum length of a basis with respect to $\{A,B,C\}$ is two, since every spanning tree intersects $\{A, B, C \}$, but the maximum length with respect to $\{a,A,B \}$ is~3, since the complement of  $\{a,A,B \}$, namely $\{b,c,C \}$ is also a basis.

The following theorem gives another characterization of binary matroids.

\begin{theorem}\label{theo:binarymat}
Let $M=(E, \cB)$ be a matroid, $F\subset E$,
$x_1, y_1\in F$ $(x_1\neq y_1)$ and $x_2, y_2\notin F$ $(x_2\neq y_2)$ such that $\cB=\Big\{F, F\Delta\{x_1, x_2\},F\Delta\{y_1, y_2\}, F\Delta\{x_1, x_2\}\Delta\{y_1, y_2\}, F\Delta\{x_1, y_2\}, F\Delta\{y_1, x_2\}\Big\}$. The matroid $M$ is not binary and no binary matroid has minor isomorphic to $M$.
\end{theorem}
\proof
Handling $x_2$ over $x_1$ gives the set 
\bea
\cB_{x_2x_1}=\cB\Delta\Big\{F\Delta\{x_1,x_2\},F\Delta\{y_1,y_2\}\Delta\{x_1,x_2\}\Big\}=\Big\{F,F\Delta\{y_1,y_2\},F\Delta\{x_1,y_2\},F\Delta\{y_1,x_2\}\Big\},
\eea
on which the exchange axiom is not satisfied. In fact, if the EA was true and $x_1\in \big(F\Delta\{y_1,x_2\}\setminus F\Delta\{x_1,y_2\}\big)$, we should find $x\in (F\Delta\{x_1,y_2\}\setminus F\Delta\{y_1,x_2\}\big)$ such that $F\Delta\{y_1,x_2\}\Delta\{x_1,x\}\in \cB_{x_2x_1}$. The possible values of $x$ are $y_1$ or $y_2$ and then $F\Delta\{y_1,x_2\}\Delta\{x_1,x\}=F\Delta\{x_1,x_2\}$ or $F\Delta\{y_1,x_2\}\Delta\{x_1,x\}=F\Delta\{y_1,x_2\}\Delta\{x_1,y_2\}$. None of them belong to $\cB_{x_2x_1}$. 
\qed

To be binary, a matroid $M$ must not have any $U_{2,4}$ minors. Theorem \ref{theo:binarymat} meets this condition because the basis set of $U_{2,4}$ is of the form $$\cB=\Big\{F, F\Delta\{x_1, x_2\},F\Delta\{y_1, y_2\}, F\Delta\{x_1, x_2\}\Delta\{y_1, y_2\}, F\Delta\{x_1, y_2\}, F\Delta\{y_1, x_2\}\Big\}.$$

\begin{lemma}\label{lemma:atleastoneortwo}
Let $M=(E, \cB)$ be a matroid defined by its base set, $F\in \cB$, $x, y, x', y'\in E$ and $x, x'\in F$. If $F\Delta\{x, y\}, F\Delta\{x', y'\}\in\cB$ and $x\neq x'$, then $y=y'$ or we have two possible  cases: $F\Delta\{x, y\}\Delta\{x', y'\}\in \cB$ or $F\Delta\{x, y'\}$ and $F\Delta\{x', y\}$ belong to $\cB$.
Furthermore if $\cB$ contains an element of the form $F\Delta\{x,y\}\Delta\{x',y'\}$ then $F\Delta\{x,y\}, F\Delta\{x',y'\}\in\cB$
or $F\Delta\{x,y'\}, F\Delta\{x',y\}\in \cB$.
\end{lemma}

\proof
Consider $F\Delta\{x, y\}, F\Delta\{x', y'\}\in\cB$, where $x\neq x'$. As
$x'\in F\Delta\{x, y\}\setminus F\Delta\{x', y'\}$ then there is $a\in F\Delta\{x', y'\}\setminus F\Delta\{x, y\}$ such that $(F\Delta\{x, y\})\Delta \{x', a\}\in \cB$. Therefore $a=x$ if $y=y'$ and then $(F\Delta\{x, y\})\Delta \{x', x\}=F\Delta\{x', y\}=F\Delta\{x', y'\}\in \cB$. In case $y\neq y'$ we have $a=x, y'$. 
This implies that $F\Delta\{x', y\}$ or $F\Delta\{x, y\}\Delta\{x', y'\}$ belong to $\cB$. Furthermore $x\in F\Delta\{x', y'\}\setminus F\Delta\{x, y\}$ implies that there is $b\in F\Delta\{x, y\}\setminus F\Delta\{x', y'\}$ such that $ F\Delta\{x', y'\}\Delta \{x, b\}\in \cB$. In the same way as above if $y=y'$ then $b=x'$ and then 
$(F\Delta\{x', y'\})\Delta \{x, x'\}=F\Delta\{x, y'\}=F\Delta\{x, y\}\in \cB$.
For $y\neq y'$ we have $b=x', y$ implying that $F\Delta\{x, y'\}$ or $F\Delta\{x, y\}\Delta\{x', y'\}$ belong to $\cB$.
This concludes the proof of the first part of this lemma. The proof for the second part is similar, as it considers $F, F\Delta\{x,y\}\Delta\{x',y'\}\in \cB$ and applies the SEA.
\qed
\begin{proposition}\label{prop:reducelength2}
Let $M=(E, \cB)$ be a binary matroid. Assume that $\cB$ is described by $F\subset E$ and any element $F'\neq F$ in $\cB$ has a maximum length of $2$. Then there is a sequence of handle slides sending $M$ to a matroid with a single element in its basis set.
\end{proposition}

\proof
If $F'$ is the only element of $\cB$ with the maximum length, any other element of $\cB$ has the form $F\Delta\{a,b\}$ with $a\in F$ and $b\notin F$. Also $F'$ is of the form $F'=F\Delta\{x_1,x_2\}\Delta\{y_1,y_2\}$ with $x_1, y_1\in F$ $(x_1\neq y_1)$ and with $x_2, y_2\notin F$ $(x_2\neq y_2)$. As $x_1\in F\setminus F'$ there is $x'\in F'\setminus F$ such that $F\Delta\{x_1,x'\}\in \cB$. Let us assume that $x'=x_2$ without losing generality. As a result, $F, F\Delta\{x_1,x_2\}, F\Delta\{x_1,x_2\}\Delta\{y_1,y_2\}\in \cB$. We can also assert that, $F\Delta\{y_1,y_2\}\in \cB$ or $F\Delta\{x_1,y_2\}, F\Delta\{y_1,x_2\}\in \cB$. Indeed, if $F\Delta\{y_1,y_2\}\notin \cB$ as $x_2\in F\Delta\{x_1,x_2\}\Delta\{y_1,y_2\}\setminus F$, the EA implies that there is $x\in F\setminus F\Delta\{x_1,x_2\}\Delta\{y_1,y_2\}$ such that $F\Delta\{x_1,x_2\}\Delta\{y_1,y_2\}\Delta\{x_2,x\}\in \cB$ and this implies that $x=y_1$. Therefore, $F\Delta\{x_1,y_2\}\in \cB$. Similarly, by taking $y_1\in F\setminus F\Delta\{x_1,x_2\}\Delta\{y_1,y_2\}$, we can prove that $F\Delta\{y_1,x_2\}\in \cB$. We obtain the following two conditions: 
the four elements $F, F\Delta\{x_1,x_2\}, F\Delta\{y_1,y_2\}, F\Delta\{x_1,x_2\}\Delta\{y_1,y_2\}$ are in $\cB$ or the five bases  $F, F\Delta\{x_1,x_2\}, F\Delta\{x_1,y_2\},F\Delta\{y_1,x_2\}, F\Delta\{x_1,x_2\}\Delta\{y_1,y_2\}$ belong to $\cB$. The above computation shows that if $\cB$ contains an element of the form $F\Delta\{x_1,x_2\}\Delta\{y_1,y_2\}$ then we get $F\Delta\{x_1,x_2\}, F\Delta\{y_1,y_2\}\in \cB$ or $F\Delta\{x_1,y_2\}, F\Delta\{y_1,x_2\}\in \cB$. 

The matroid $M_{x_1x_2}$ will lose the bases $F\Delta\{x_1,x_2\}$ and
$F\Delta\{x_1,x_2\}\Delta\{y_1,y_2\}$ but may gain base of the form $F\Delta\{x_1,x_2\}\Delta\{a,b\}$; $a\in F$ and $b\notin F$. This is possible only when $F\Delta\{a,b\}\in \cB$ and $a\neq x_1$, $b\neq x_2$. From our previous result $F\Delta\{a,x_2\}, F\Delta\{x_1,b\}\in \cB$ and taking $(M_{x_2x_1})_{x_2,a}$ or $(M_{x_2x_1})_{b,x_1}$ will kill the new base $F\Delta\{x_1,x_2\}\Delta\{a,b\}$. Repeating the same computation after a sequence of handle slides will yield a binary matroid $M'=(E, \cB')$ in which every element except from $F$ is of the form $F\Delta\{p,q\}$ with $p\in F$ and $q\notin F$ and can be killed in $M'_{qp}$ and if any other element $F\Delta\{r,s\}\Delta\{p,q\}$ is created, it will still be removed by handle slides applying the process described earlier.

In case $M$ has another element $F'=F\Delta\{a_1,b_1\}\Delta\{b_1,b_2\}$ of length $2$, such that $F\Delta\{a_1,b_1\}$, $F\Delta\{b_1,b_2\}\in \cB$ then we can possibly get in $M_{x_2x_1}$ an element of length $3$ that is of the form $F\Delta\{a_1,b_1\}\Delta\{a_2,b_2\}\Delta\{x_1,x_2\}$ that will be removed as well by computing $(M_{x_2x_1})_{a_2a_1}$ because the base set $F\Delta\{b_1,b_2\}\Delta\{x_1,x_2\}\in \cB_{x_2x_1}$. If one of the elements $F\Delta\{a_1,b_1\}$ or $F\Delta\{b_1,b_2\}\notin \cB$ then applying the previous results, 
$F\Delta\{a_1,b_2\}, F\Delta\{b_1,a_2\}\in \cB$. Since $F\Delta\{a_1,b_1\}\Delta\{a_2,b_2\}=F\Delta\{a_1,b_2\}\Delta\{a_2,b_1\}$ we can proceed as above. In conclusion any new element of length $3$ that we could create can be removed through a sequence of handle slides. The same process can be extended to any element with a length greater than $3$ through induction.

Finally a sequence of handle slides allows us to remove all the bases elements of $M$ except one.
\qed

\begin{theorem}\label{theo:binarymatroid}
Let $M=(E, \cB)$ be a binary matroid. There is a sequence of handle slides that sends $M$ to a binary matroid with a single element in the basis set.
\end{theorem}

\proof
Let us proceed by induction on the maximum length $n$ of the elements in $\cB$.
The Proposition \ref{prop:reducelength2} solves the cases $n=1,2$. Assume that $\cB$ contains an element $F'$ of length $3$ of the form 
$F\Delta\{a,b\}\Delta\{p,q\}\Delta\{r,s\}$. As a result, we can claim with confidence that $F\Delta\{a,b\}\in \cB$
and $F\Delta\{p,q\}\in \cB$.
In fact, because $a\in F\setminus F'$, there is $x\in F'\setminus F$ such that $F\Delta\{a,x\}\in \cB$ and we can set $x=b$. Furthermore, $p\in F\Delta\{a,b\}\setminus F'$ and thus there is $y\in F'\setminus F\Delta\{a,b\}$ such that $F\Delta\{a,b\}\Delta\{p,y\}\cB$ and we can take $y=q$. Besides, $F\Delta\{p,x\}\in \cB$ for some $x=b,q,s$ as $p\in F\setminus F'$. If $x=p,s$, we are done. Otherwise, $F\Delta\{p,b\}\in \cB$. In addition, $F\Delta\{a,b\}\Delta\{p,q\}\Delta\{b,x\}\in \cB$ for some $x\in F\setminus F\Delta\{a,b\}\Delta\{p,q\}$. As a result, $x=a,p$. We are done if $x=a$, otherwise $F\Delta\{a,q\}\in \cB$. Then we have $F\Delta\{a,q\}, F\Delta\{p,b\}\in \cB$, which also solve the problem because we can interchange the role between $p$ and $a$.
Letting $F\Delta\{a,b\}$ play the role of $F$ we can apply the second result in Lemma  \ref{lemma:atleastoneortwo} and obtain $F\Delta\{a,b\}\Delta\{p,q\}, F\Delta\{a,b\}\Delta\{r,s\}\in \cB$ or $F\Delta\{a,b\}\Delta\{p,s\}, F\Delta\{a,b\}\Delta\{r,q\}\in \cB$. Continuing in the same manner as before, we don't need to study the two cases since both of them will lead to the same conclusion. So we get $F\Delta\{a,b\}, F\Delta\{p,q\} F\Delta\{a,b\}\Delta\{p,q\}, F\Delta\{a,b\}\Delta\{r,s\}\in \cB$ and $F'$ is killed in $M_{pq}$ or $M_{rs}$. Finally, it is straightforward to demonstrate that we can send $M$ to a binary matroid with a single element in its basis set using a series of handle slides.
\qed

The binary matroid given by this theorem is known as the canonical binary matroid, and the single element of the basis set corresponds to the set of edges of a spanning tree of the graph $G$ if $M=M(G)$ is a connected graphic matroid.

\section{Canonical binary $\Delta$-matroids}
The goal of this section is to solve the earlier introduced conjecture. As a result, we investigate and address the main result of this paper. 
\begin{proposition}\label{prop:canonicalint}
Let $D=(E,\cF)$ be a $\Delta$-matroid. There is a sequence of handle slides that sends $D$ to a $\Delta$-matroid whose feasible set contains only one element of minimum size and one element of maximum size.
\end{proposition}
\proof
Let $D_l=(E, \cF_l)$ and
$D_u=(E, \cF_u)$ be respectively the lower and the upper matroids associated with a binary $\Delta$-matroid $D$. It is obvious to see that if in $D_{ab}$ ($a, b\in E$; $a\neq b$) we create a new feasible $F'$, then $\rk (D_l)\leq|F'|\leq\rk(D_u)$. This inequality implies that a sequence of handle slides performed on $D$ will always result in a $\Delta$-matroid $D'$ of the same rank with $D$ and $\rk (D_l)=\rk (D'_l)$, $\rk (D_u)=\rk (D'_u)$. 

There is a sequence of handle slides from Theorem \ref{theo:binarymatroid} that sends $D_l$ and $D_u$ to canonical binary matroids $C_l$ and $C_u$ respectively. The result is obtained by applying to $D$ the sequence of handle slides sending $D u$ to $C_u$, followed by the sequence sending $D_l$ to $C_l$. In fact, the first sequence of handle slides results to a $\Delta$-matroid, in which feasible set contains a single  element $F$ with maximum size. Furthermore, Theorem \ref{the:spanind} 
asserts that the lower matroid's basis is independent in the upper matroid, and that every basis in $D_l$ is contained in $F$. Applying the handle slide sequence from $D_l$ to $C_l$ will end the proof because any pair $a,b$ in this sequence belongs to $F$ and will prevent us from creating new feasibles of size $|F|$.
 
\qed
\begin{theorem}\label{theo:deltamatcanoni}
For each binary $\Delta$-matroid $D$, there is a sequence of handle slides taking $D$ to some $D_{i,j,k,l}$, where $i$ is the size of the ground set minus the size of a largest feasible set, $l$ is the size of a smallest feasible set, $2j + k$ is difference in the sizes of a largest and a smallest feasible
set. Furthermore, the handle slides can take $D$ to $D_{i,j,0,l}$ if $D$ is even and to $D_{i,0,k,l}$ if $D$ is odd.
\end{theorem}
\proof
Applying Proposition \ref{prop:canonicalint} to a binary $\Delta$-matroid $D$ yields a $\Delta$-matroid $D'$ with a single feasible $F$ that is included in every other feasible of $D'$. The inclusion is due to Theorem \ref{the:spanind} in \cite{AS2021}.
Therefore, $D$ is isomorphic to the $\Delta$-matroid obtained by taking the direct sum of the $\Delta$-matroid $D'/F\oplus (F,\{F\})$. With $D'/F$ the $\Delta$-matroid resulting from the contraction in $D'$ all the elements of $F$ and $(F,\{F\})$ the matroid whose ground set is $F$ and $F$ is the only element of its basis. Because the binary $\Delta$-matroid $D'/F$ has the empty set as a feasible, and therefore from Theorem \ref{theo:binarydeltaempty} it can be sent to a canonical binary $\Delta$-matroid $D_{i,j,k}$. In addition, $(F,\{F\})=\oplus_{e\in F}(\{e\}, \{\{e\}\})$ and setting $|F|=l$, we can clearly send (by handle slides) $D$ to $D_{i,j,k,l}$, the $\Delta$-matroid
consisting of the direct sum of $i$ copies of $({e}, \{\emptyset\})$, $j$ copies of $(\{e, f\}, \{\emptyset, \{e, f\}\})$, $k$ copies of $(\{e\}, \{\emptyset, \{e\}\})$ and $l$ copies of
$(\{e\}, \{\{e\}\})$.

Furthermore, $D$ is even if and only if $D'$ is even, i.e.
$k=0$, and otherwise by a sequence of handle slides we can send $D$
to $D_{i,0,k,l}$ with $k>0$.
\qed

\begin{figure}[htb]
\begin{center}
\unitlength=1.85mm
       \begin{picture}(16.25,16.25)
           \put(0.00,0.00){\line(1,0){10.00}}
           \put(0.00,0.00){\line(0,1){10.00}}
           \put(0.00,0.00){\line(1,1){10.00}}
           \put(10.00,0.00){\line(0,1){10.00}}
           \put(0.00,0.00){\circle*{2.50}}
           \put (0.00,-3.00) {$2$}
           \put(10.00,0.00){\circle*{2.50}}
           \put (10.00,-3.00) {$3$}
           \put(0.00,10.00){\circle*{2.50}}
           \put (0.00,12.00) {$1$}
           \put(10.00,10.00){\circle*{2.50}}
           \put (10.00,12.00) {$4$}

\end{picture}
\end{center}
   \caption{The graph of adjacency matrix $A$.}
   \label{fig:adj}
\end{figure}
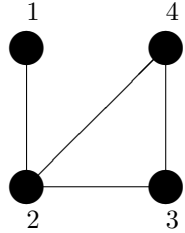

\begin{example}
The pair $D=(E, \cF)$, with $E=\{1,2,3,4\}$ and $$\cF=\big\{\{1\}, \{2\}, \{1,2,3\}, \{1,2,4\}, \{1,3,4\}, \{2,3,4\}\big\},$$ is a binary $\Delta$-matroid. In fact taking a twist with $\{1\}$ gives the $\Delta$-matroid in which the feasible sets are the invertible
submatrices of the adjacency matrix
$A$ of the graph in Figure \ref{fig:adj}. The vertex $i$ of the graph corresponds to  row/column $i$ of $A$.

$$
  A =\begin{pmatrix}
0&1&0&0\\
1&0&1&1\\0&1&0&1\\0&1&1&0
\end{pmatrix}.
$$
We now consider the following handle slides:
$$\mathcal{F}_{12} = \{ \{2\}, \{1,2,3\}, \{1,2,4\}, \{2,3,4\}\},$$
$$(\mathcal{F}_{12})_{34} = \{\{2\}, \{1,2,4\}, \{2,3,4\}\},$$
$$\big((\mathcal{F}_{12})_{34}\big)_{13}= \{ \{2\}, \{2,3,4\} \},$$
and therefore by a sequence of handle slides we send $D$ to $(\{1\},\{\emptyset\})\oplus(\{3,4\},\{\emptyset, \{3,4\}\})\oplus(\{2\},\{\{2\}\})$ which is isomorphic to the canonical $\Delta$-matroid $D_{1,1,0,1}$.
\end{example}

\begin{figure}
    \centering
    \includegraphics{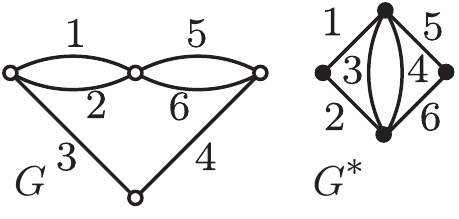}
    \caption{$G$ and $G^*$ on the sphere.}
    \label{fig:sphere}
\end{figure}

\begin{figure}
    \centering
    \includegraphics{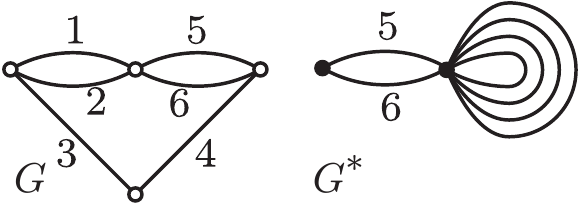}
    \caption{$G$ and $G^*$ on torus or projective plane.}
    \label{fig:torus}
\end{figure}

Cellularly embedded graphs on surfaces provide a rich source for binary $\Delta$-matroids~\cite{MR4193372}.
The lower matroid of such a $\Delta$-matroid is the cycle matroid of the embedded graph $G$, while the upper matroid is isomorphic to the co-cycle matroid of $G^*$, the geometric dual of $G$.
For example consider the graph $G$ in Figures~\ref{fig:sphere},\ref{fig:torus}. If $G$ is embedded on the sphere, then the cycle matroid of $G$ and the co-cycle matroid of $G^*$ are isomorphic and the upper and lower matroid are identical and the feasible sets $\cF$ of the corresponding $\Delta$-matroid are the edge sets of the spanning trees of $G$, 
\bea
\cF_1 = &&\Big\{\{1,3,4\}, \{1,3,5\}, \{1,3,6\}, \{1,4,5\}, \{1,4,6\},
 \{2,3,4\}, \{2,3,5\}, \{2,3,6\}, \{2,4,5\},
\{2,4,6\},\cr && \{3,4,5\}, \{3,4,6\} \Big\}.
\eea

However, if $G$ is embedded on a torus such that its geometric dual becomes $G^*$ as in Figure~\ref{fig:torus}, then we get in addition to $\cF_1$ the two sets
$$\{1,2,3,4,5\}, \{1,2,3,4,6\} \},$$
so the corresponding $\Delta$-matroid has 
$\cF_2 = \cF_1 \cup \{1,2,3,4,5\}, \{1,2,3,4,6\} \} $ as feasibles.
The same $G$ and $G^*$ could be interpreted as coming from a map on the projective plane, in which case the set
$\{1,2,3,4\}$ also becomes feasible and 
$\cF_3 = \cF_2 \cup \{1,2,3,4\}$ is again the collection of feasible sets of a  $\Delta$-matroid, this time the feasible sets are of even and odd parity.

Let us examine the normal forms obtainable from $\cF_1, \cF_2,\cF_3$ via a sequence of handle slides. For $\cF_1$ we can arrive at the edge set of a spanning tree of $G$ of our choice, say $\{1,3,4\}$.  
Our canonical $\Delta$-matroid becomes
 $D_{3,0,0,3} =(\{2\},\{\emptyset\})\oplus(\{5\},\{\emptyset\})\oplus(\{6\},\{\emptyset\})\oplus(\{1\},\{\{1\}\})\oplus(\{3\},\{\{3\}\})\oplus(\{4\},\{\{4\}\})$. 
 
 A normal form for a $\Delta$-matroid with
 $\cF_2$ as feasibles is
 $D_{1,1,0,3} =(\{6\},\{\emptyset\})\oplus(\{4,5\},\{\emptyset, \{4,5\}\})\oplus(\{1\},\{\{1\}\})\oplus(\{3\},\{\{3\}\})\oplus(\{4\},\{\{4\}\})$. 
 
 For $\cF_3$, our example derived from a non-orientable map, we get
 $D_{1,0,2,3} =(\{6\},\{\emptyset\})\oplus(\{2\},\{\emptyset, \{2\}\})
 \oplus(\{5\},\{\emptyset, \{5\}\})
 \oplus(\{1\},\{\{1\}\})\oplus(\{3\},\{\{3\}\})\oplus(\{4\},\{\{4\}\})$.

\bibliographystyle{alpha}
\bibliography{references}

\vspace{0.5cm}

\end{document}